\documentclass[11pt]{amsart}

\usepackage{amsfonts}
\usepackage{amsmath}
\usepackage{amssymb}
\usepackage{amsthm}

\usepackage{graphicx}

\usepackage{color}

\usepackage{hyperref}

\usepackage{paralist}

\usepackage[margin=1.08in]{geometry}

\newcommand\set[2]{\left\{ {#1} \ : \ {#2} \right\}}

\newcommand\R{\mathbb{R}}               
\newcommand\Z{\mathbb{Z}}               

\newcommand\id{\mathrm{Id}}

\renewcommand\int{\mathrm{int}\,}

\newcommand\rk{\mathrm{rk}}

\newcommand\hyp{\Lambda}

\newcommand\x{\mathbf{x}}
\newcommand\e{\mathbf{e}}
\newcommand\1{\mathbf{1}}
\renewcommand\id{\mathrm{Id}}
\newcommand\q{\mathbf{q}}

\newtheorem{thm}{Theorem}[section]

\theoremstyle{definition}

\newtheorem{example}[thm]{Example}

\title[On derivative cones of polyhedral cones]{On the derivative cones of
polyhedral cones}
\author{Raman Sanyal}
\address{Raman Sanyal, Department of Mathematics, FU Berlin, Arnimallee 2,
14195 Berlin, Germany}
\email{sanyal@math.fu-berlin.de}

\keywords{polyhedral cones, hyperbolic polynomials, derivative cones,
spectrahedra, elementary symmetric functions, polymatroids}
\subjclass[2000]{90C22, 52A41, 52B99}

\date{\today}
\thanks{ This research was supported by a Miller Research Fellowship at
UC Berkeley.}

\begin{document}

\begin{abstract}
    Hyperbolic polynomials elegantly encode a rich class of convex cones that
    includes polyhedral and spectrahedral cones. Hyperbolic polynomials are
    closed under taking polars and the corresponding cones -- the derivative
    cones -- yield relaxations for the associated optimization problem and
    exhibit interesting facial properties.  While it is unknown if every
    hyperbolicity cone is a section of the positive semidefinite cone, it is
    natural to ask whether spectrahedral cones are closed under taking
    derivative cones.  In this note we give an affirmative answer for
    polyhedral cones by exhibiting an explicit spectrahedral representation
    for the first derivative cone.  We also prove that higher polars do not
    have an determinantal representation which shows that the problem for
    general spectrahedral cones is considerably more difficult.
\end{abstract}

\maketitle


\section{Introduction}

A homogeneous polynomial $p(\x) \in \R[x_1,\dots,x_n]$ is \emph{hyperbolic}
with respect to a direction $\e \in \R^n$ if $p(\e) > 0$ and for every $\q \in \R^n$
the univariate polynomial
\[
    \lambda \ \mapsto \ p(\q + \lambda \e )
\] 
has only real roots. The study of hyperbolic polynomials originated in the theory
of (hyperbolic) differential equations  where it was noted by
G\r{a}rding~\cite{garding59} that the connected component of $\R^n \setminus V(p)$
containing $\e$ is an (open) convex cone. Its closure is called the
\emph{hyperbolicity cone} $\hyp_{p,\e}$ associated to the pair $(p,\e)$.  G\"uler
\cite{guler} introduced hyperbolic polynomials to convex optimization and showed
that interior-point methods can be effectively employed for hyperbolicity cones.
More precisely, he showed that $-\log p(\x)$ serves as a self-concordant barrier
function for $\hyp_{p,\e}$.

Two of the  most important subclasses of hyperbolicity cones are the
\emph{polyhedral} and the \emph{spectrahedral} cones. For a full-dimensional
polyhedral cone 
\[
    P \ = \ \{ \x \in \R^n \;:\; \ell_i(\x) \ge 0 \text{ for all } i =
                1,\dots,d \},
\]
presented as the intersection of finitely many halfspaces given by linear
forms $\ell_1(\x),\dots,\ell_d(\x)$ it is clear that the polynomial
\[
    p(\x) \ = \ \ell_1(\x) \cdot \ell_2(\x) \cdots \ell_d(\x)
\]
is hyperbolic with respect to $\e \in P$ in the interior of $P$ and $\hyp_{p,\e} =
P$. A spectrahedral cone $S$ is the intersection of the cone of positive
semidefinite matrices with a linear subspace~\cite{rg}. By
choosing a basis for the linear subspace, the cone can be presented by
\begin{equation}\label{eqn:det}
    S \ = \ \{ \x \in \R^n \;:\; A(\x) = x_1 A_1 + x_2 A_2 + \cdots + x_n A_n \
    \succeq 0 \}
\end{equation}
where $A_1,\dots,A_n \in \R^{d \times d}$ are symmetric matrices and $A(\x)
\succeq 0$ denotes positive semidefiniteness. In particular, $S$ has a presentation
such that $A(\e)$ is positive definite for all interior points $\e \in S$.
The corresponding hyperbolic polynomial is given by
\begin{equation}\label{eqn:detform}
    p(\x) \ = \ \det A(\x) \ = \ \det\, x_1 A_1 + x_2 A_2 + \cdots + x_n A_n
\end{equation}
and it is easy to see that $p(\q - \lambda \e)$ is the characteristic polynomial
of a symmetric matrix and hence has only real zeros. These two examples in
particular show the richness of the class of convex cones given by hyperbolic
polynomials. However, it is still unknown whether this class of convex cones
is strictly larger than spectrahedral cones.  In a stronger form Peter
Lax~\cite{lax} conjectured that every ternary hyperbolic polynomial has a
determinantal representation in the sense of \eqref{eqn:detform} with $A(\e)$
positive definite. This conjecture was answered in the positive by Lewis et
al.~\cite{lrp05} using results of Helton and Vinnikov~\cite{hv}. In higher
dimensions the Lax conjecture is generally false as can be inferred from a
count of parameters: hyperbolic polynomials constitute a set with non-empty
interior inside the space of forms of fixed degree~\cite{nuij68}. In contrast,
symmetric matrices of order $d$ yield determinantal forms of type
\eqref{eqn:detform} that are constrained to a proper algebraic subset of the
space of homogeneous polynomials of degree $d$.

Spectrahedral cones constitute a very nice class of cones from a convex
geometric point of view. For instance, spectrahedral cones are \emph{facially
exposed}. This holds true for general hyperbolicity cones as was shown by
Renegar~\cite{renegar06}. In his investigation of the boundary structure, he
considered cones obtained from polarizations of $p(\x)$ which yield
relaxations of $\hyp_{p,\e}$ while sharing parts of the facial structure with
the original cone.  For a homogeneous polynomial $p(\x)$, the \emph{$i$-th
polar} with respect to $\e$ is 
\[ 
    R^i_\e p(\x) \ := \ \sum_{|\alpha| = i} \e^\alpha \frac{\partial^i}{\partial
    \x^\alpha} p(\x),
\] 
that is, the degree $i$ part in a Taylor expansion of $p(\x)$ around $\e$. In
particular, the first polar is given by $R^1_\e p(\x) = \sum_i e_i
\frac{\partial}{\partial x_i}p(\x) = \langle \nabla p(\x), \e \rangle$ and it
follows from Rolle's theorem that $R^i_\e p$ is hyperbolic with respect to
$\e$ whenever $p(\x)$ is. The corresponding cones are the \emph{derivative
cones} $\hyp^i_{p,\e} = \hyp_{R^i_\e p,\e}$ of $(p,\e)$.

Towards a deeper understanding of the geometry of hyperbolicity cones, it is
natural to ask if the subclass of spectrahedral cones is closed under taking
derivative cones. The main result of this note is that this is true for the
first derivative cone of a polyhedral cone. 
\begin{thm} \label{thm:main1}
    Let $P = \set{ \x \in \R^n }{ \ell_i(\x) \ge 0 \text{ for } i \in [d] }$ be a
    full-dimensional polyhedral cone. Let $\e \in \int P$ and assume that
    $\ell_i(\e) = 1$ for all $i$. Then the first derivative cone is given by all
    $\x \in \R^n$ satisfying
    \[
        \begin{bmatrix}
            \ell_1(\x) + \ell_n(\x)  & \ell_n(\x) & \cdots & \ell_n(\x) \\
            \ell_n(\x) &  \ell_2(\x) + \ell_n(\x)  &  \cdots & \ell_n(\x) \\
            \vdots  & \vdots & \ddots & \vdots \\
            \ell_n(\x)  & \ell_n(\x) & \cdots & \ell_{n-1}(\x) + \ell_n(\x) \\
            \end{bmatrix} \ \succeq \ 0.
        \]
\end{thm}

Let us illustrate the theorem with an example that motivated this research.

\begin{example}
    The $3$-dimensional \emph{halfcube} $H_3$ is the subpolytope of the $3$-cube
    $[-1,1]^3$ induced by the vertices with an even number of $-1$ coordinates. It
    is the affine slice of the polyhedral cone below with $\{t = 1\}$. The
    corresponding polynomial $p$ is hyperbolic with respect to $\e = (1,0,0,0)$.
    The same affine section of its first derivative cone is the $3$-dimensional
    \emph{elliptope}~\cite{lp} (also known as the \emph{Samosa}).
    \begin{center}
        \begin{tabular}{c@{\hspace{1cm}}c}
        \includegraphics[width=4cm]{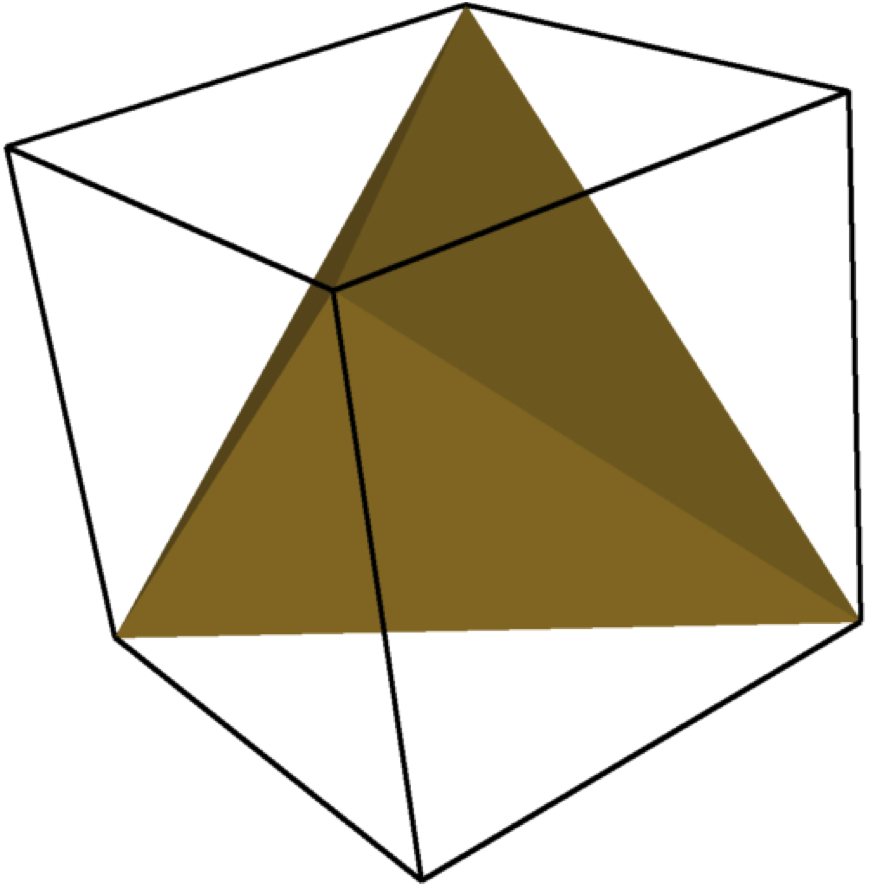} & 
        \includegraphics[width=4cm]{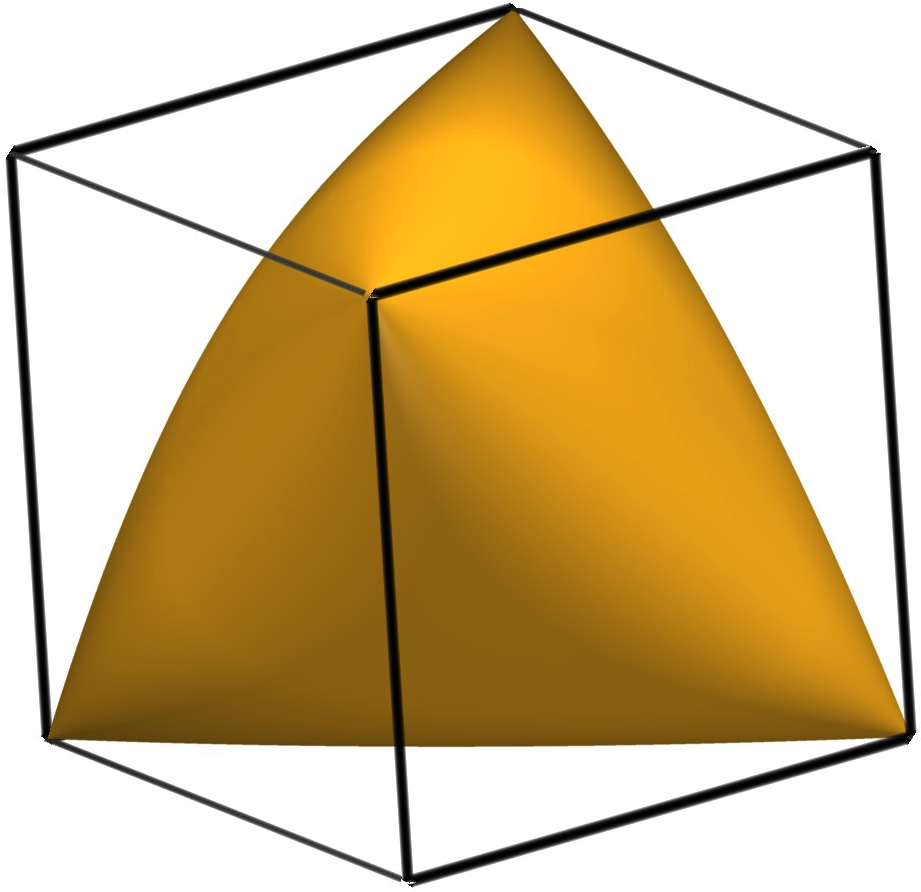} \\
        $
        \begin{array}{r@{\ \  \le \ \ }l}
            x  \ - \ y \ - \ z  & t\\
            -x \ + \ y \ - \ z & t\\
            -x \ - \ y \ + \ z & t\\
            x  \ + \ y \ + \ z  & t\\
        \end{array}
        $
        &
        $
        \begin{bmatrix}
            2t-2x & t-x-y-z & t-x-y-z \\
            t-x-y-z & 2t-2y & t-x-y-z \\
            t-x-y-z & t-x-y-z & 2t-2z \\
        \end{bmatrix} \ \succeq \ 0
        $
    \end{tabular}
    \end{center}
    The Samosa has $6$ one-dimensional faces corresponding to the edges of the
    halfcube. The first polar is the equation \emph{Cayley cubic surface} $t^3
    - tx^2 - ty^2 - tz^2 + 2xyz$. Unfortunately, this nice connection between
    halfcubes and elliptopes does not hold in higher dimensions. \hfill
    $\diamond$
\end{example}

In fact, the above theorem asserts the much stronger statement that the first
polar of a product of linear forms satisfies the generalized Lax conjecture, that
is, the first polar of $p(\x) = \prod_i \ell_i(\x)$ equals the determinant of the
given matrix.  However, this is the only polar for which this is true.
\begin{thm} \label{thm:main2}
    Derivative cones of spectrahedra do not satisfy the generalized Lax
    conjecture.
\end{thm}

The prototypical example of a polyhedral cone is given by $\R^d_{\ge 0}$ with
hyperbolic polynomial $p(\x) = E_d(\x) = x_1x_2\cdots x_d$, the $d$-th elementary
symmetric function. In particular, $E_d(\x)$ is hyperbolic with respect to $\e =
\1$ and $R^i_\1 E_d(\x) = E_{d-i}(\x)$ are elementary symmetric functions. The
relation of the corresponding cones to semidefinite programming was considered in
Zinchenko~\cite{zinchenko} who showed that $\hyp^i_{p,\e}$ is a
\emph{spectrahedral shadow}, i.e.~the projection of a spectrahedral cone. However,
this is a weaker statement in that spectrahedral cones are not closed under
projection and images under projection have much weaker properties such as
non-exposed faces. The main focus of~\cite{zinchenko} was a description for the
dual cone to $\hyp_{E_{n-1},\1}$ which can be inferred from
Theorem~\ref{thm:main1}. The remainder of the paper is concerned with the proof of
the two theorems above.

{\bf Acknowledgments.} The author would like to thank the organizers of the
AIM workshop ``Convex algebraic geometry, optimization and applications'' where this question was raised. We would like to
thank Thilo R\"orig and Philipp Rostalski for fruitful discussions and Bernd
Sturmfels for support and encouragement.

\section{Proofs via realizations of uniform matroids}

Before we come to the proof of Theorem~\ref{thm:main1}, let us reduce the
statement to the case of the positive orthant $\R^d_{\ge 0}$ with hyperbolic
polynomial $E_d(\x) = x_1x_2 \cdots x_d$.  To that end, note that for $p(\x) =
\prod_i^d \ell_i(\x)$, the $i$-th polar satisfies
\[
R^i_\e p(\x) \ = \ i! p(\e) E_{d-i}(\ell_1(\x),\dots,\ell_d(\x))
\]
where $E_{d-i}(y_1,\dots,y_d)$ is the $d-i$-th elementary symmetric function.
Geometrically, this means that the $i$-th derivative cone of $P = \R^d_{\ge0}
\cap L$, where $L$ is a linear subspace, equals the intersection of $L$ with the $i$-th derivative cone of
$\R^d_{\ge 0}$. The individual linear forms $\ell_i(\x)$ are called the
\emph{eigenvalues} of hyperbolic polynomial $p(\x)$ and the above claim holds
true for general hyperbolic polynomials~\cite[Prop.~18]{renegar06}.

While Theorem~\ref{thm:main1} can be proved inductively by verifying that the
determinant of the shown matrix equals the first polar of $p(\x)$ with respect to $\1$,
we will make use of a beautiful connection between hyperbolic polynomials
and polymatroids due to Gurvits~\cite{gurvits05} and spectrahedra and
realizations of polymatroids due to Choe et al.~\cite{choe04}. Recall that a
\emph{polymatroid} on the ground set $[n] = \{1, 2, \dots, n\}$ is a monotone
submodular function $\rk : 2^{[n]} \rightarrow \Z_{\ge 0}$. That is, for any
$I, J \subseteq [n]$ we have $0 \le \rk(I) \le \rk(I \cup J)$ and
\[
    \rk(I \cup J) \ \le \ \rk(I) + \rk(J) - \rk( I \cap J ).
\]
If, in addition, $\rk(I) \le |I|$ for all $I \subseteq [n]$, then $\rk$
defines a \emph{matroid}; see~\cite[Sect.~12.1]{oxley} for more information.
Matroids and polymatroids are combinatorial abstractions of subspace
arrangements: Let $V_1, V_2,\dots, V_n$ be linear subspaces of a given vector
space $V$, then 
\[
    I \ \mapsto \dim \sum_{i \in I} V_i
\]
yields a polymatroid. The case of matroids corresponds to line arrangements.
Polymatroids from subspace arrangements are called \emph{realizable}
polymatroids and it is well known that not all polymatroids are
realizable~\cite[Ch.~6]{oxley}. 
\begin{example}[Uniform matroids]
    For every $0 \le k \le n$, the function $\rk_k(I) = \min(k, |I|)$ defines
    a matroid, called the \emph{uniform matroid} $U_{k,n}$. A realization of
    $U_{k,n}$ is given by any generic choice of $n$ lines in a $k$-dimensional
    vector space.
\end{example}

The connection between hyperbolic polynomials and polymatroids was first observed
by Gurvits~\cite{gurvits05}. For a univariate polynomial $p(t) \in \R[t]$, let
$\mathrm{mult}_0(p(t))$ be the order of vanishing at the origin. For a subset $I
\subseteq [n]$, denote by $\chi_I \in \{0,1\}^n$ its characteristic vector.
\begin{thm}[{\cite[Fact~5.3]{gurvits05}}] \label{thm:hyp_polym}
    Let $p(\x)$ be a $d$-form hyperbolic with respect to $\e$ and assume that
    $\R^n_{\ge 0} \subseteq \hyp_{p,\e}$. Then the map $\rk_p: 2^{[n]} \rightarrow
    \Z$ given by
    \[
        I \ \mapsto d - \mathrm{mult}_0(p(\chi_I + \lambda \e))
    \]
    defines a polymatroid.
\end{thm}
        
It is easy to see that for $e = \1$, the polymatroid associated to the $k$-th
elementary symmetric function $E_k(\x)$ is exactly the uniform matroid
$U_{k,n}$. 

The connection to spectrahedral cones is that the form \eqref{eqn:det} yields a
realization of the associated polymatroid. To see this, let us assume that the
matrix map
\[
    A(\x) \ = \ x_1 A_1 + x_2 A_2 + \cdots + x_n A_n
\]
with symmetric $A_1, \dots, A_n \in \R^{d \times d}$ satisfies
$A(\e) = \id$. If $\det A(\x)$ satisfies the condition of
Theorem~\ref{thm:hyp_polym}, then every $A_i$ is positive semidefinite and
hence is of the form $A_i = L_iL_i^T$ for some real matrix $L_i \in \R^{d
\times k_i}$. Now, $p(\chi_I + \lambda \e) = \det(A(\chi_I) + \lambda \id)$ is
the characteristic polynomial of $A(\chi_I)$. The order of vanishing is the
nullity of $L_I^T$, where $L_I = (L_i : i \in I)$ is the concatenation of
$L_i$ for $i \in I$.  Hence $d - \det(A(\chi_I) - \lambda \id)$ is the rank of
$L_I$. The observation of Choe et al.~\cite{choe04} is that the hyperbolic
polynomial $p(\x) = \det(A(\x))$ contains information about the realization
given by $L = (L_1,L_2,\dots,L_n)$. In our context, it is sufficient to assume
that $L_i$ is actually a single vector. Let $X = \mathrm{diag}(x_1,\dots,x_n)$, then
\[
    p(\x) \ = \ \det A(\x) \ = \ \det( L \cdot X \cdot L^T )  \ = \ 
    \sum\set{ \det(L_I)^2 \x^I }{ I \subseteq [n], |I| = d }
\]
where the last equality follows from the Binet-Cauchy theorem.

We are now able to give the proofs of Theorems \ref{thm:main1} and
\ref{thm:main2}: For the former, we need to find a realization of the uniform
matroid $U_{n-1,n}$ in the form of a matrix $L \in \R^{(n-1) \times n}$ such
that every maximal minor is $\pm 1$. Such a \emph{unimodular} realization for
$U_{n-1,n}$ is given by $L = (\id, -\1)$.  Theorem~\ref{thm:main2} follows
from the well-known fact that $U_{k,n}$ has no unimodular realization unless
$k \le 1$ or $k \ge n-1$.  Indeed, if $2 \le k \le n-2$ and $L \in \R^{k
\times n}$ matrix giving an unimodular realization of $U_{k,n}$ then, without
loss of generality, $L = (\id,L^\prime)$ with $L^\prime \in \R^{k \times
(n-k)}$. The conditions on the maximal minors of $L$ forces $L^\prime$ to have
all minors $\pm 1$ which is only possible iff $L^\prime$ has at most one row
or column. This proves Theorems~\ref{thm:main1} and \ref{thm:main2}.

Let us remark that this does \emph{not} exclude the possibility that the
hyperbolicity cones for $E_k(\x)$ are slices of the PSD cone. Indeed, it would
be sufficient to find an polynomial $g(\x)$ hyperbolic with respect to $\e =
\1$ such that $\hyp_{E_k,\e} \subseteq \hyp_{g,\e}$ and $gE_k$ has a symmetric
determinantal representation  with $A(\e) \succ 0$. We close by showing this
for $E_2(\x)$. To this end, consider the spectrahedral cone given by 
\[
A(\x) \ = \
\begin{bmatrix} 
    E_1(\x) & x_1 & x_2 & \cdots & x_{n-1} \\
     x_1 & E_1(\x) \\
     x_2 & & E_1(\x) \\
     \vdots &&& \ddots \\
     x_{n-1}  &&&& E_1(\x) \\
\end{bmatrix} \ \succeq \ 0.
\]
It is easy to see that $A(\1)$ is positive definite and
\[
\det A(\x) \ = 2 E_1(\x)^{n-2} E_2(\x).
\]
The arrangement of $2$-planes is not very well understood. As of now, we do
not know spectrahedral representations for all elementary symmetric functions
but we conjecture that they exist.
    

\bibliographystyle{siam}
\bibliography{RenegarOfPolytope}

\begin{thebibliography}{10}

\bibitem{choe04}
{\sc Y.-B. Choe, J.~G. Oxley, A.~D. Sokal, and D.~G. Wagner}, {\em Homogeneous
  multivariate polynomials with the half-plane property}, Adv. in Appl. Math.,
  32 (2004), pp.~88--187.
\newblock Special issue on the Tutte polynomial.

\bibitem{garding59}
{\sc L.~G\r{a}rding}, {\em An inequality for hyperbolic polynomials}, J. Math.
  Mech., 8 (1959), pp.~957--965.

\bibitem{guler}
{\sc O.~G{\"u}ler}, {\em Hyperbolic polynomials and interior point methods for
  convex programming}, Math. Oper. Res., 22 (1997), pp.~350--377.

\bibitem{gurvits05}
{\sc L.~Gurvits}, {\em Combinatorial and algorithmic aspects of hyperbolic
  polynomials}.
\newblock preprint, \url{math.CO/0404474}, April 2005.

\bibitem{hv}
{\sc J.~W. Helton and V.~Vinnikov}, {\em Linear matrix inequality
  representation of sets}, Comm. Pure Appl. Math., 60 (2007), pp.~654--674.

\bibitem{lp}
{\sc M.~Laurent and S.~Poljak}, {\em On the facial structure of the set of
  correlation matrices}, SIAM J. Matrix Anal. Appl., 17 (1996), pp.~530--547.

\bibitem{lax}
{\sc P.~D. Lax}, {\em Differential equations, difference equations and matrix
  theory}, Comm. Pure Appl. Math., 11 (1958), pp.~175--194.

\bibitem{lrp05}
{\sc A.~S. Lewis, P.~A. Parrilo, and M.~V. Ramana}, {\em The {L}ax conjecture
  is true}, Proc. Amer. Math. Soc., 133 (2005), pp.~2495--2499 (electronic).

\bibitem{nuij68}
{\sc W.~Nuij}, {\em A note on hyperbolic polynomials}, Math. Scand., 23 (1968),
  pp.~69--72 (1969).

\bibitem{oxley}
{\sc J.~G. Oxley}, {\em Matroid theory}, Oxford Science Publications, The
  Clarendon Press Oxford University Press, New York, 1992.

\bibitem{rg}
{\sc M.~Ramana and A.~J. Goldman}, {\em Some geometric results in semidefinite
  programming}, J. Global Optim., 7 (1995), pp.~33--50.

\bibitem{renegar06}
{\sc J.~Renegar}, {\em Hyperbolic programs, and their derivative relaxations},
  Found. Comput. Math., 6 (2006), pp.~59--79.

\bibitem{zinchenko}
{\sc Y.~Zinchenko}, {\em On hyperbolicity cones associated with elementary
  symmetric polynomials}, Optim. Lett., 2 (2008), pp.~389--402.

\end{thebibliography}

\end{document}